\documentclass[final,1p,times,number]{elsarticle}

\usepackage{amsmath}
\usepackage{amssymb}
\usepackage{amsfonts}
\usepackage{graphicx}
\usepackage{epstopdf}
\usepackage{color}
\usepackage{bm}
\usepackage{multirow}
\usepackage{natbib}
\usepackage{hyperref}
\usepackage{cleveref}

\usepackage{color}
\usepackage{float}
\usepackage{caption}

\usepackage[caption=false]{subfig}

\usepackage{ifpdf}
\ifpdf%
\usepackage{pdflscape}
\else
\usepackage{lscape}
\fi

\hypersetup{
	colorlinks,
	linkcolor={red!50!black},
	citecolor={blue!50!black},
	urlcolor={blue!80!black}
}

\def\a{\mbox{\boldmath $a$}}

\def\f{\mbox{\boldmath $f$}}

\def\g{\mbox{\boldmath $g$}}

\def\m{\mbox{\boldmath $m$}}

\def\x{\mbox{\boldmath $x$}}

\def\0{\mbox{\boldmath $0$}}

\biboptions{sort&compress}

\begin{document}

\begin{frontmatter}

\title{Semi-implicit Structure Preserving Method for The Landau-Lifshitz Equation}

\author[XJTLU]{Changjian Xie\corref{cor}}
\cortext[cor]{Corresponding author.} 
\ead{Changjian.Xie@xjtlu.edu.cn}



\address[XJTLU]{School of Mathematics and Physics, Xi'an-Jiaotong-Liverpool University, Re'ai Rd. 111, Suzhou, 215123, Jiangsu, China.}

\begin{abstract}
A critical challenge inherent to the projection method applied to the Landau-Lifshitz equation is the deficiency of rigorous theoretical justifications for the stability of its projection step. To mitigate this limitation, we introduce a semi-implicit numerical scheme, which is formulated on the basis of the first-order Backward Differentiation Formula (BDF) incorporated with one-sided extrapolation and a Crank-Nicolson-type norm-preserving procedure. This proposed scheme exhibits three fundamental characteristics: structure preservation, numerical stability, and first-order accuracy in time. In practical implementations, the scheme not only ensures stable computation and adheres to the norm constraint but also guarantees the uniqueness of the numerical solution, thereby providing substantial facilitation for the theoretical analysis of the normalizing step.

\end{abstract}


\end{frontmatter}

\section{Introduction}

Ferromagnetic materials are widely employed in data storage technologies, primarily owing to the bistable characteristics of their intrinsic magnetic ordering, commonly known as magnetization. The temporal evolution of magnetization in such systems is fundamentally governed by the Landau-Lifshitz-Gilbert (LLG) equation \cite{Landau1935On,Gilbert:1955}. This equation encompasses two essential components that regulate magnetization dynamics: a gyromagnetic term responsible for energy conservation, and a damping term that accounts for energy dissipation.
The LLG equation forms a vector-valued nonlinear system subject to a pointwise constant-magnitude constraint on the magnetization vector. Extensive research has been devoted to developing efficient and numerically stable algorithms for micromagnetic simulations, as documented in review articles including \cite{kruzik2006recent,cimrak2007survey}. Among the various numerical approaches, semi-implicit schemes have attracted substantial attention due to their favorable stability properties and avoidance of sophisticated nonlinear solvers \cite{alouges2006convergence, gao2014optimal, Xie2018}. For instance, our group previously proposed a second-order backward differentiation formula (BDF2) scheme based on one-sided interpolation \cite{Xie2018}, which requires solving a three-dimensional linear system with non-constant coefficients at each time step. A rigorous convergence analysis confirming the second-order accuracy of this BDF2 scheme was presented in \cite{jingrun2019analysis}. Alternatively, Alouges et al. \cite{alouges2006convergence} introduced a linearly implicit scheme using the tangent space to enforce the magnetization constraint, albeit with only first-order temporal accuracy. More recently, Lubich et al. \cite{Lubich2021} constructed and analyzed high-order BDF schemes for the LLG equation. Unconditional unique solvability for such semi-implicit schemes has been rigorously established in \cite{jingrun2019analysis,Lubich2021}; nevertheless, their convergence analysis typically imposes a restriction that the time step is proportional to the spatial mesh size. Despite these advances, the practical implementation and assessment of high-order numerical methods in micromagnetic modeling remain active research directions.
Although numerous studies have investigated structure-preserving algorithms for the LLG equation, see \cite{wang2001gauss}, \cite{qing_2023}, \cite{xiaoli_2024}, \cite{JEONG2010613}, and \cite{An2016OptimalEE}, such methods are not directly derived from the governing equation itself. Instead, starting from the requirement of norm preservation, they typically introduce nonlinear projection steps, adopt Lagrange multipliers, or embed the constraint indirectly into the discrete formulation.
A critical drawback remains prevalent in existing approaches:
most available methods rely on projection techniques to enforce the norm-preserving constraint. While such projection preserves the orientation of the magnetization, it normalizes its magnitude. Due to its inherent nonlinearity, this projection step introduces substantial difficulties in the theoretical analysis of stability and convergence.
To overcome this limitation, this work presents a first-order accurate numerical scheme for the LLG equation, combining a semi-implicit BDF1 method for stable evolution and a Crank–Nicolson-type predictor to preserve the normalization constraint naturally. In this way, no artificial constraint is imposed; instead, the structure of the model itself is exploited to achieve stable and consistent discrete solutions. We also perform extensive numerical experiments to validate the numerical stability of the proposed scheme, which further demonstrate its accuracy and stability properties.

The rest of this paper is structured as follows. \cref{sec: numerical scheme} first reviews the micromagnetic model, then provides a detailed formulation of the proposed numerical scheme, which modifies the conventional BDF1 projection method into a BDF1 Crank–Nicolson-type structure-preserving algorithm. \cref{sec:experiments} reports comprehensive numerical results, including verification of temporal and spatial accuracy in both one-dimensional (1D) and three-dimensional (3D) settings.
Finally, concluding remarks and perspectives for future work are given in \cref{sec:conclusions}.

\section{The LLG equation}
\label{sec: numerical scheme}


Micromagnetics, as a fundamental discipline investigating the magnetic behaviors of ferromagnetic materials at the mesoscopic scale, relies heavily on the Landau-Lifshitz-Gilbert (LLG) equation as its core governing framework. Proposed to describe the spatiotemporal dynamics of magnetization in ferromagnetic systems, this equation rigorously integrates two indispensable physical mechanisms that dominate magnetization evolution: gyromagnetic precession and dissipative relaxation \cite{Landau1935On,Brown1963micromagnetics}. Its nondimensionalized form, which facilitates analytical derivation and numerical implementation, is given by
\begin{align}\label{c1-large}
{\m}_t =-{\m}\times{\bm h}_{\text{eff}}-\alpha{\m}\times({\m}\times{\bm h}_{\text{eff}}),
\end{align}
where $\partial_t \m$ is the time derivative for the magnetization and ${\bm h}_{\text{eff}}$ is the term from the energy variation.

This equation is subject to the homogeneous Neumann boundary condition

\begin{equation}\label{boundary-large}
\frac{\partial{\m}}{\partial {\bm \nu}}\Big|_{\partial \Omega}=0,
\end{equation}
where \(\Omega \subset \mathbb{R}^d\) (\(d=1,2,3\)) denotes the bounded spatial domain occupied by the ferromagnetic material, and \(\bm \nu\) represents the unit outward normal vector on the domain boundary \(\partial \Omega\). It is noteworthy that this boundary condition is physically consistent for isolated ferromagnetic systems, as it inherently ensures the absence of magnetic surface charge—an essential prerequisite for accurately modeling unperturbed magnetic dynamics.

To fully comprehend the LLG equation, it is critical to clarify the physical nature of its key components. The magnetization field \(\m: \Omega \to \mathbb{R}^3\) is a three-dimensional vector field constrained by the pointwise condition \(|\m|=1\), a fundamental property rooted in the quantum mechanical alignment of electron spins within ferromagnetic materials. Regarding the right-hand side of \cref{c1-large}, the first term characterizes the gyromagnetic precession phenomenon, wherein magnetic moments undergo precessional motion around the effective magnetic field \(\bm h_{\text{eff}}\). The second term accounts for dissipative relaxation, with the parameter \(\alpha > 0\) denoting the dimensionless Gilbert damping coefficient that quantifies the rate of energy transfer from the magnetic subsystem to the lattice structure.

The effective magnetic field \(\bm h_{\text{eff}}\), a pivotal quantity in the LLG equation, is mathematically derived as the functional derivative of the Gibbs free energy functional \(F[\m]\) with respect to the magnetization field, i.e., \(\bm h_{\text{eff}} = -\delta F/\delta \m\). For the purpose of simplifying subsequent analytical derivations and numerical discretizations without compromising the core physical essence, we adopt the Gibbs free energy functional in the following simplified form:

\begin{equation}\label{LL-Energy}
F[\m] = \int_\Omega |\nabla\m|^2 \mathrm{d}\x . 
\end{equation}

Under this formulation, the effective field \(\bm h_{\text{eff}}\) can be explicitly derived by evaluating the functional derivative, leading to the concise expression \(\bm h_{\text{eff}} = \Delta\m\)—a result that captures the core physical contributions to the effective field within the scope of the present study. Substituting this explicit form of \(\bm h_{\text{eff}}\) into the general LLG equation (\cref{c1-large}), we obtain a simplified yet equivalent formulation of the governing equation:

\begin{align}\label{eq-5}
\m_t=-\m\times\Delta\m-\alpha\m\times(\m\times\Delta\m).
\end{align}

For the sake of facilitating stable numerical discretization, \cref{eq-5} can be further simplified by leveraging vector algebra identities and the inherent constraint of the magnetization field. Specifically, utilizing the vector triple product identity \(\a \times ({\bm b} \times {\bm c}) = (\a \cdot {\bm c}){\bm b} - (\a \cdot {\bm b}){\bm c}\) and the pointwise constraint \(|\m| = 1\) (from which \(\m \cdot \partial_t \m = 0\) can be directly deduced via time differentiation), we derive an alternative equivalent form of the LLG equation:

\begin{equation}\label{eq-model}
\m_t=\alpha  \Delta\m-\m\times\Delta\m+\alpha |\nabla \m|^2\m.
\end{equation}

A critical distinction in existing research on LLG equation discretization lies in the choice of formulation adopted. Specifically, numerous prior studies (e.g., \cite{qing_2023}) have focused on the simplified form given by \cref{eq-model}, treating this formulation as the heat flow of harmonic maps and incorporating \(|\nabla \m|\) as a Lagrange multiplier to enforce the magnetization constraint. In stark contrast to these conventional approaches, the present work establishes its numerical framework explicitly based on the LLG formulation in \cref{eq-5}. This deliberate choice not only avoids the indirect constraint-handling strategies associated with \cref{eq-model} but also directly retains the intrinsic physical structure of the original LLG equation—an essential foundation that underpins the novelty and rigor of the numerical scheme proposed in this study.
 
\section{Proposed method}

For the construction of the structure-preserving method, we set \(\alpha=0\) in \cref{eq-5}, yielding the simplified LLG equation without damping:
\begin{align}\label{eq-alpha-0}
\m_t=-\m\times\Delta\m.
\end{align}

To lay the theoretical foundation for the proposed scheme, we first consider the simple linear vectorial equation

\begin{align}\label{eq-CN}
\m_t=-\m\times \a,
\end{align}

where \(\a^T=(a_1,a_2,a_3)\) denotes a constant vector. Applying the Crank-Nicolson method to \cref{eq-CN} leads to the discrete formulation:

\begin{align}\label{eq-CN_1}
	\frac{\m_h^{n+1}-\m_h^n}{\Delta t}=-\frac{\m_h^{n+1}+\m_h^n}{2} \times \a,
\end{align}

This discrete scheme is norm-preserving, which can be verified by taking the inner product of both sides with \(\m_h^{n+1}+\m_h^n\), resulting in:

\begin{align*}
	\|\m_h^{n+1}\|_2=\|\m_h^n\|_2.
\end{align*}

Rearranging \cref{eq-CN_1} into a compact linear system form, we obtain:

\begin{align*}
	\begin{pmatrix}
	1&\frac12 \Delta t a_3&-\frac12 \Delta t a_2\\
	-\frac12 \Delta ta_3&1&\frac12 \Delta ta_1\\
	\frac12 \Delta ta_2&-\frac12 \Delta t a_1&1
	\end{pmatrix}\begin{pmatrix}
	m_1^{n+1}\\
	m_2^{n+1}\\
	m_3^{n+1}
	\end{pmatrix}=\begin{pmatrix}
	m_1^n+\frac12 \Delta t(a_2 m_3^n-a_3 m_2^n)\\
	m_2^n+\frac12 \Delta t(a_3 m_1^n-a_1 m_3^n)\\
	m_3^n+\frac12 \Delta t(a_1m_2^n-a_2 m_1^n)
	\end{pmatrix}
\end{align*}

This can be further rewritten as a matrix-vector multiplication of the form:

\begin{align*}
	\begin{pmatrix}
	m_1^{n+1}\\
	m_2^{n+1}\\
	m_3^{n+1}
	\end{pmatrix}=\begin{pmatrix}
	1&\frac12 \Delta t a_3&-\frac12 \Delta t a_2\\
	-\frac12 \Delta ta_3&1&\frac12 \Delta ta_1\\
	\frac12 \Delta ta_2&-\frac12 \Delta t a_1&1
	\end{pmatrix}^{-1}\begin{pmatrix}
	1&-\frac12 \Delta t a_3 & \frac12 \Delta t a_2\\
	\frac12 \Delta t a_3&1& -\frac12 \Delta t a_1\\
	-\frac12 \Delta t a_2&\frac12 \Delta t a_1&1
	\end{pmatrix}\begin{pmatrix}
	m_1^n\\
	m_2^n\\
	m_3^n
	\end{pmatrix}=A\begin{pmatrix}
	m_1^n\\
	m_2^n\\
	m_3^n
	\end{pmatrix}
\end{align*}

where the matrix \(A\) is explicitly given by:

\begin{align*}
	A=\frac{1}{S}\begin{pmatrix}
	1+\beta^{2}a_1^{2}&-2\beta a_3+\beta^{2}a_1a_2&2\beta a_2+\beta^{2}a_1a_3\\
	2\beta a_3+\beta^{2}a_1a_2&1+\beta^{2}a_2^{2}&-2\beta a_1+\beta^{2}a_2a_3\\
	-2\beta a_2+\beta^{2}a_1a_3&2\beta a_1+\beta^{2}a_2a_3&1+\beta^{2}a_3^{2}
	\end{pmatrix}
\end{align*}

Here, \(S=\det(A)=1+\beta^2(a_1^2+a_2^2+a_3^2)\) and \(\beta=\frac{\Delta t}{2}\). It can be rigorously proven that \(\|A\|_2=1\), confirming the stability of the Crank-Nicolson scheme for \cref{eq-CN}.

Building on this theoretical basis, we propose three structure-preserving schemes for \cref{eq-alpha-0} as follows:

\begin{itemize}
	\item Scheme I: Explicit treatment of \(\Delta \m\),
	\begin{align}\label{eq-CN_2}
	\frac{\m_h^{n+1}-\m_h^n}{\Delta t}=-\frac{\m_h^{n+1}+\m_h^n}{2} \times \Delta_h \m_h^n,
	\end{align}
	This scheme is subject to a CFL-type condition for numerical stability.
		\item Scheme II: Implicit treatment of \(\Delta \m\),
	\begin{align}\label{eq-CN_3}
	\frac{\m_h^{n+1}-\m_h^n}{\Delta t}=-\frac{\m_h^{n+1}+\m_h^n}{2} \times \Delta_h \m_h^{n+1},
	\end{align}
	This formulation introduces computational challenges due to the high complexity of the resulting nonlinear system.
	\item Scheme III: Semi-implicit method,
	\begin{align}\label{eq-CN_4}
	\frac{\m_h^{n+1}-\m_h^n}{\Delta t}=-\frac{\m_h^{n+1}+\m_h^n}{2} \times \Delta_h \g_h,
	\end{align}
	where \(\g_h^{s}=(I-\Delta t \Delta_h)^{-1}\m_h^s\) for \(s=n,n+1\). Specifically, this scheme can be expressed as the linear system:
	\begin{align*}
	\begin{pmatrix}
	1&\frac12 \Delta t \Delta_hg_3&-\frac12 \Delta t \Delta_hg_2\\
	-\frac12 \Delta t\Delta_hg_3&1&\frac12 \Delta t\Delta_hg_1\\
	\frac12 \Delta t\Delta_hg_2&-\frac12 \Delta t \Delta_hg_1&1
	\end{pmatrix}\begin{pmatrix}
	m_1^{n+1}\\
	m_2^{n+1}\\
	m_3^{n+1}
	\end{pmatrix}=\begin{pmatrix}
	m_1^n+\frac12 \Delta t(\Delta_hg_2 m_3^n-\Delta_hg_3 m_2^n)\\
	m_2^n+\frac12 \Delta t(\Delta_hg_3 m_1^n-\Delta_hg_1 m_3^n)\\
	m_3^n+\frac12 \Delta t(\Delta_hg_1m_2^n-\Delta_hg_2 m_1^n)
	\end{pmatrix}
	\end{align*}
	The key design trick lies in the evaluation time of \(\g_h\) (i.e., at \(t_n\) or \(t_{n+1}\)). When \(\g_h\) is evaluated at \(t_n\), the scheme becomes:
	\begin{align*}
	\begin{pmatrix}
	1&\frac12 \Delta t \Delta_hg_3^n&-\frac12 \Delta t \Delta_hg_2^n\\
	-\frac12 \Delta t\Delta_hg_3^n&1&\frac12 \Delta t\Delta_hg_1^n\\
	\frac12 \Delta t\Delta_hg_2^n&-\frac12 \Delta t \Delta_hg_1^n&1
	\end{pmatrix}\begin{pmatrix}
	m_1^{n+1}\\
	m_2^{n+1}\\
	m_3^{n+1}
	\end{pmatrix}=\begin{pmatrix}
	m_1^n+\frac12 \Delta t(\Delta_hg_2^n m_3^n-\Delta_hg_3^n m_2^n)\\
	m_2^n+\frac12 \Delta t(\Delta_hg_3^n m_1^n-\Delta_hg_1^n m_3^n)\\
	m_3^n+\frac12 \Delta t(\Delta_hg_1^n m_2^n-\Delta_hg_2^n m_1^n)
	\end{pmatrix}
	\end{align*}
	This variant is numerically verified to outperform the CFL constraint of Scheme I. Based on the above analysis, we finally propose the following structure-preserving method for the original LLG equation:
    \begin{equation}\label{eq-proposed-scheme}
\left\{ 
\begin{aligned}
&\frac{{\tilde{\m}}_h^{n+1} -  {\m}_h^n}{k}
= -\m^{n}_h\times \Delta_h \tilde{{\m}}_h^{n+1}-\alpha \m_h^n \times (\m^{n}\times \Delta_h \tilde{{\m}}_h^{n+1}),\\
&\frac{{{\m}}_h^{n+1} -  {\m}_h^n}{k}
=  -\left(\frac{\m_h^{n+1}+\m_h^n}{2}\right)\times \Delta_h \tilde{\m}_h^{n+1}-\alpha\left(\frac{\m_h^{n+1}+\m_h^n}{2}\right)\times (\m_h^n \times \Delta_h \tilde{\m}_h^{n+1}).
\end{aligned}
\right.
\end{equation} 
    \end{itemize}
    The first step to get the stable $\tilde{\m}_h^{n+1}$, but without the norm preservation. The second step due to cross product with the average of $\m_h^{n+1}$ and $\m_h^n$, such that it can preserve the norm constraints.

\section{Numerical experiments}
\label{sec:experiments}

\subsection{Accuracy tests}

To quantitatively evaluate the numerical accuracy of the proposed scheme, analytical exact solutions are derived for both one-dimensional (1D) and three-dimensional (3D) scenarios, serving as benchmark solutions for error quantification of the method applied to \cref{eq-5}.
For the 1D case, the exact magnetization solution \(\m_e\) is given by:

\begin{equation*}
\m_e=\left(\cos(\cos(\pi x))\sin t, \sin(\cos(\pi x))\sin t, \cos t\right)^T,
\end{equation*}

whereas the corresponding 3D exact solution is formulated as:

\begin{equation*}
\m_e=\left(\cos(XYZ)\sin t, \sin(XYZ)\sin t, \cos t\right)^T,
\end{equation*}

with \(X=x^2(1-x)^2\), \(Y=y^2(1-y)^2\), and \(Z=z^2(1-z)^2\) denoting the spatial modulation terms in the 3D domain.

These exact solutions are strictly consistent with the governing equation \cref{eq-5} when the forcing term is defined as \(\f_e=\partial_t \m_e+\m_e \times \Delta \m_e+\alpha\m_e\times (\m_e \times \Delta \m_e)\). Additionally, they satisfy the homogeneous Neumann boundary condition specified in \cref{boundary-large}, ensuring full consistency with the simulation constraints imposed in the present study.

To isolate the temporal approximation error from the effects of spatial discretization, the spatial resolution in the 1D accuracy test is fixed at \(h=5\times 10^{-4}\)—a sufficiently fine grid that renders spatial discretization error negligible compared to temporal approximation error. The Gilbert damping parameter is set to \(\alpha=0.01\), and all simulations are conducted until the final time \(T=0.1\). Under this configuration, the measured numerical error primarily reflects the inaccuracy introduced by temporal discretization of the proposed scheme.

The 3D temporal accuracy test is inherently constrained by spatial resolution, as excessively fine spatial grids lead to prohibitive computational costs. To balance the contributions of spatial and temporal errors and ensure reliable accuracy assessment, a coordinated refinement strategy is adopted for the spatial mesh sizes (\(h_x, h_y, h_z\)) and temporal step-size (\(\Delta t\)), tailored to the first-order convergence of the proposed method: \(\Delta t=h_x^2=h_y^2=h_z^2=h^2=T/N_0\). Here, \(N_0\) denotes a refinement level parameter, whose specific values are provided in the subsequent numerical results section. Consistent with the 1D test, the Gilbert damping parameter is set to \(\alpha=0.01\), and the final simulation time \(T\) is specified in the corresponding result tables. The first-order temporal accuracy of the proposed scheme is verified by the numerical results presented in \cref{tab-a-v-3} (for the 1D test) and \cref{tab-a-10-time-3D-Q-2} (for the 3D test), respectively.

\begin{table}[htbp]
	\centering
	\caption{The temporal accuracy in 1D test for proposed method with damping $\alpha=0.01$when $h = 5D-4$, $T=1d-1$.}\label{tab-a-v-3}
	\begin{tabular}{|c|c|c|c|}
		\hline
		$k$ & $\|\m_h-\m_e\|_\infty$ & $\|\m_h-\m_e\|_2$ &$\|\m_h-\m_e\|_{H^1}$\\
		\hline
		2d-2 &0.018930003955949&0.011905019571153&0.055826921067779 \\
		1d-2 &0.010028587881722&0.006173875738341&0.029296196267458 \\
		5d-3 &0.005196870030573&0.003165482579077&0.015278145278317 \\
		2.5d-3 &0.002655553006303&0.001609938361519&0.007888811720037 \\
		1.25d-3 &0.001342412485260&8.133526108798887e-04&0.004027197400045 \\
        6.25d-4 &6.742979183864614e-04&4.090379381003823e-04&0.002038034696297 \\
        3.125d-4 &3.377273200392897e-04&2.051458547423506e-04&0.001025746673531 \\
		\hline 
		order &0.970286000783516&0.977447889884335&0.961185387141062\\
		\hline
	\end{tabular}
\end{table}

\begin{table}[htbp]
	\centering
	\caption{The spatial accuracy in 1D test for proposed method with damping $\alpha=0.01$ when $k= 1D-6$, $T=1d-1$.}\label{tab-a-10-space-Q-1}
	\begin{tabular}{|c|c|c|c|}
		\hline
		$h$ & $\|\m_h-\m_e\|_\infty$ & $\|\m_h-\m_e\|_2$ &$\|\m_h-\m_e\|_{H^1}$ \\
		\hline
		1/16 &4.212674325335744e-04&2.887722168727909e-04&0.002212111086235 \\
		1/24 &1.872260365181275e-04&1.277863702263473e-04&9.790353208801370e-04 \\
		1/32 &1.049626868055292e-04&7.163958605254648e-05&5.502778761147910e-04 \\
		1/48 &4.608916242324762e-05&3.159939889548999e-05&2.449630533418659e-04 \\
		1/64 &2.545990677620125e-05&1.760034075909782e-05&1.382981214094384e-04 \\
		\hline 
		order &2.023480961445598&2.017528405593854&1.999719394821664 \\
		\hline
	\end{tabular}
\end{table}

Following the temporal accuracy evaluation, spatial accuracy tests were conducted to quantify the spatial discretization performance of the proposed method. To prevent temporal errors from interfering with the assessment of spatial accuracy, the temporal step size was fixed at a sufficiently small value of \(k=10^{-6}\) for the 1D spatial accuracy test presented in \cref{tab-a-10-space-Q-1}—a choice that renders temporal discretization errors negligible compared to spatial discretization errors. For the 3D spatial accuracy test, we adopt the relationship \(k=h^2\); this configuration ensures that second-order spatial accuracy is observed, which is consistent with the convergence behavior obtained from the 1D spatial accuracy test.

\begin{table}[htbp]
	\centering
	\caption{The temporal accuracy and spatial accuracy for the proposed method with damping $\alpha=0.01$, $T=0.1$ with $k=h^2$ in 3D. }\label{tab-a-10-time-3D-Q-2}
	\begin{tabular}{|c|c|c|c|c|}
		\hline
		$k$&$h$ & $\|\m_h-\m_e\|_\infty$ & $\|\m_h-\m_e\|_2$ &$\|\m_h-\m_e\|_{H^1}$ \\
	\hline
	T/10&1/10&4.997325319509027e-04&2.885786030289455e-04&3.172841357302907e-04 \\
	T/40 &1/20&1.257686203441910e-04&7.230181520301087e-05&1.106111804599846e-04 \\
	T/57 &1/24&8.861344600608057e-05&5.087773586103405e-05&9.141039466366754e-05 \\
	T/78&1/28&6.508006096783703e-05&3.736395858299103e-05&7.960500518126427e-05 \\
	T/102 &1/32&5.005337568275703e-05&2.879318679898088e-05&7.230037732663106e-05 \\
	\hline 
	order &&0.991507777634215&0.993698915740268&0.653283961155192\\
		\hline
	&	order &1.977866420534538&1.982242514017785&1.303278322578679\\
		\hline
	\end{tabular}
\end{table}

\subsection{Norm preserving tests}
Specifically, the initial conditions employed for the 1D and 3D spatial accuracy tests are chosen as follows: for the 1D case, the initial magnetization \(\m_0\) is given by
\begin{align*}
    \m_0=\left(\cos(\cos(\pi x))\sin (0.01), \sin(\cos(\pi x))\sin (0.01), \cos (0.01)\right)^T,
\end{align*}

while for the 3D case, the initial magnetization \(\m_0\) is formulated as

\begin{align*}
\m_0=\left(\cos(XYZ)\sin (0.01), \sin(XYZ)\sin (0.01), \cos (0.01)\right)^T,
\end{align*}

where \(X=x^2(1-x)^2\), \(Y=y^2(1-y)^2\), and \(Z=z^2(1-z)^2\) are the same spatial modulation terms as defined in the 3D exact solution. The numerical results of the 1D and 3D spatial accuracy tests are presented in \cref{tab-8} and \cref{tab-9}, respectively.

\begin{table}[htbp]
	\centering
	\caption{The proposed method with damping $\alpha=0.01$ when $h = 5D-4$, $T=1d-1$ for the norm preserving test in 1D.}\label{tab-8}
	\begin{tabular}{|c|c|}
		\hline
		$k$ & $\|\|\m_h\|_2-1\|_\infty$ \\
		\hline
	    2.0D-2 & 1.554312234475219e-15 \\
		1.0D-2 &1.776356839400250e-15  \\
		5.0D-3 & 3.108624468950438e-15\\
		2.5D-3 & 4.662936703425657e-15 \\
		1.25D-3 & 6.439293542825908e-15 \\
		6.25D-4 & 7.993605777301127e-15 \\
		3.125D-4 &1.287858708565182e-14 \\
		\hline
	\end{tabular}
\end{table}

\begin{table}[htbp]
	\centering
	\caption{The norm preserving test for the proposed method with damping $\alpha=0.01$, $T=0.1$ with $k=h^2$ in 3D. }\label{tab-9}
	\begin{tabular}{|c|c|c|}
		\hline
		$k$&$h$ & $\|\|\m_h\|_2-1\|_\infty$  \\
	\hline
	T/10&1/10&1.459055098962381e-12\\
	T/40 &1/20&7.391864897954292e-13\\
	T/57 &1/24& 6.215028491851626e-13\\
	T/78&1/28& 5.242473122279989e-13\\
	\hline 
	\end{tabular}
\end{table}

\subsection{Initial conditions and stability tests}

To further verify the numerical consistency, stability, and robustness of the proposed method, different initial conditions are specified for each test case, as detailed below.
In the 1D scenario, the initial conditions are specified as follows:

\begin{align*}
    \m_0&=\left(\cos(\cos(\pi x))\sin (0), \sin(\cos(\pi x))\sin (0), \cos (0)\right)^T\\
    \m_0&=\left(\cos(\cos(\pi x))\sin (0.01), \sin(\cos(\pi x))\sin (0.01), \cos (0.01)\right)^T.
\end{align*}

\begin{figure}[htbp]
    \centering
    \subfloat[$T0=0$, $m_1$]{\includegraphics[width=0.5\linewidth]{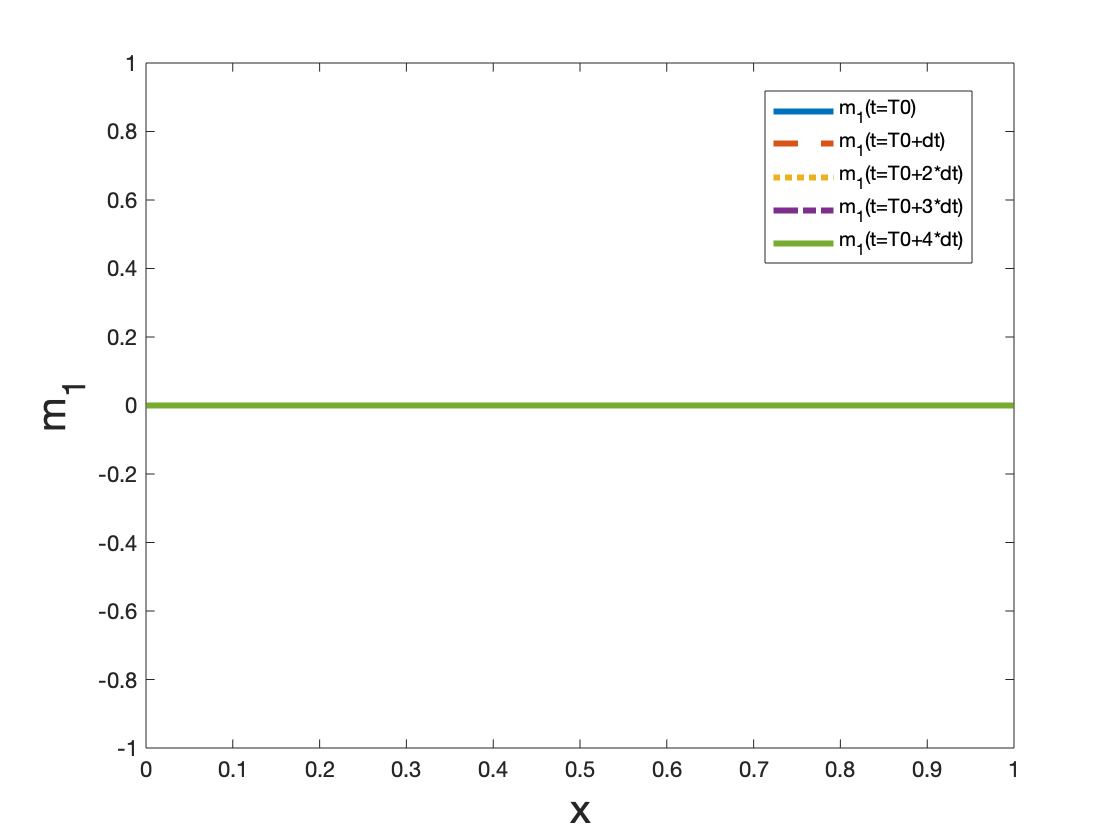}}
    \subfloat[$T0=0.01$, $m_1$]{\includegraphics[width=0.5\linewidth]{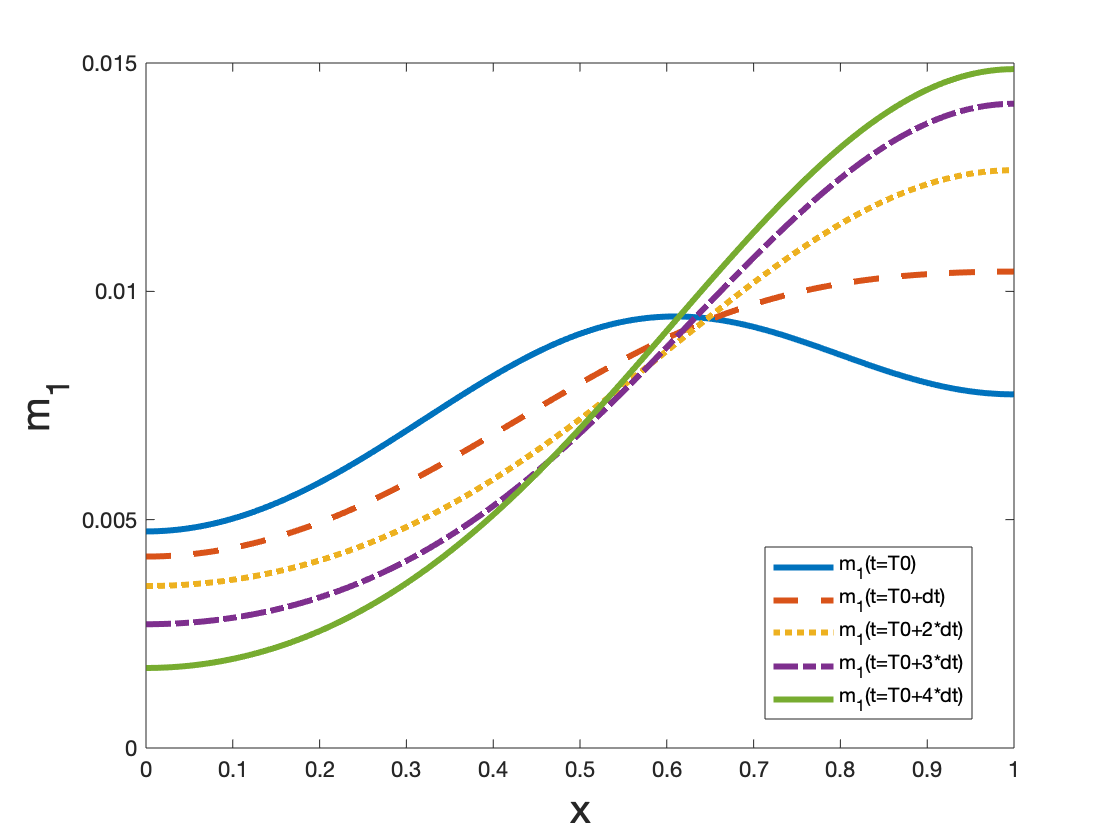}}
    \hspace{0.1in}
    \subfloat[$T0=0$, $m_2$]{\includegraphics[width=0.5\linewidth]{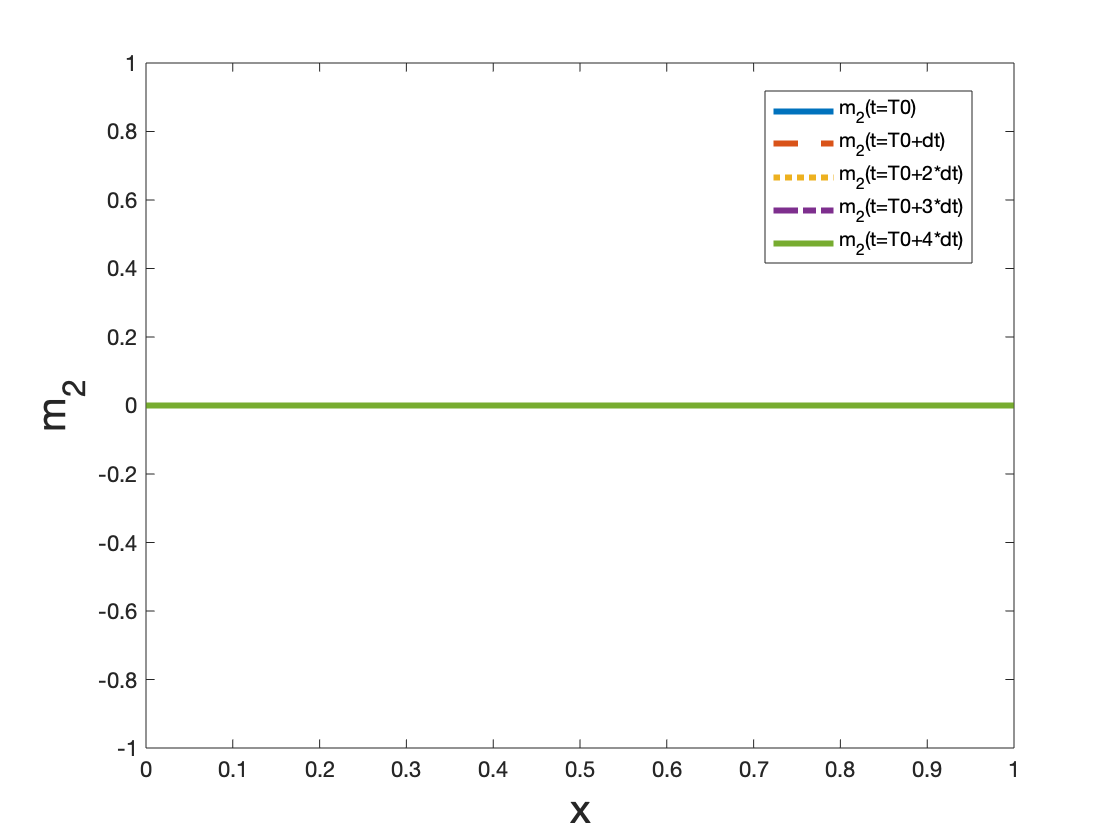}}
    \subfloat[$T0=0.01$, $m_2$]{\includegraphics[width=0.5\linewidth]{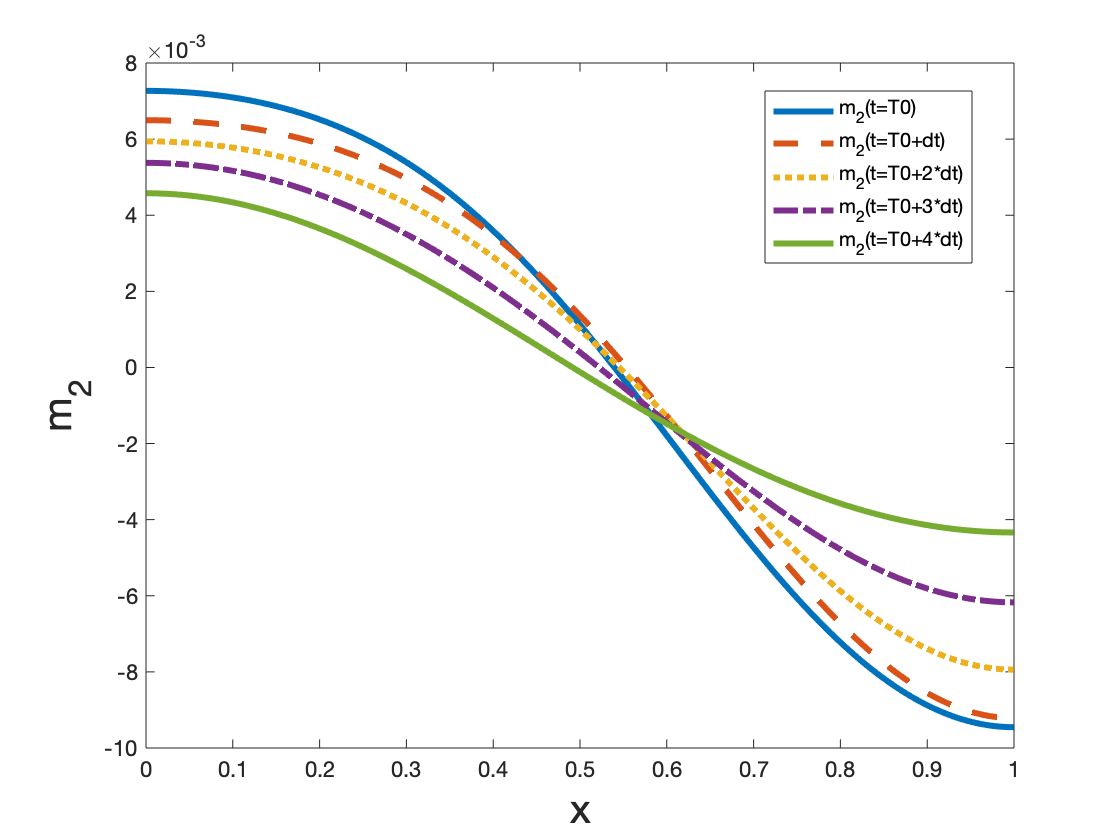}}
    \hspace{0.1in}
    \subfloat[$T0=0$, $m_3$]{\includegraphics[width=0.5\linewidth]{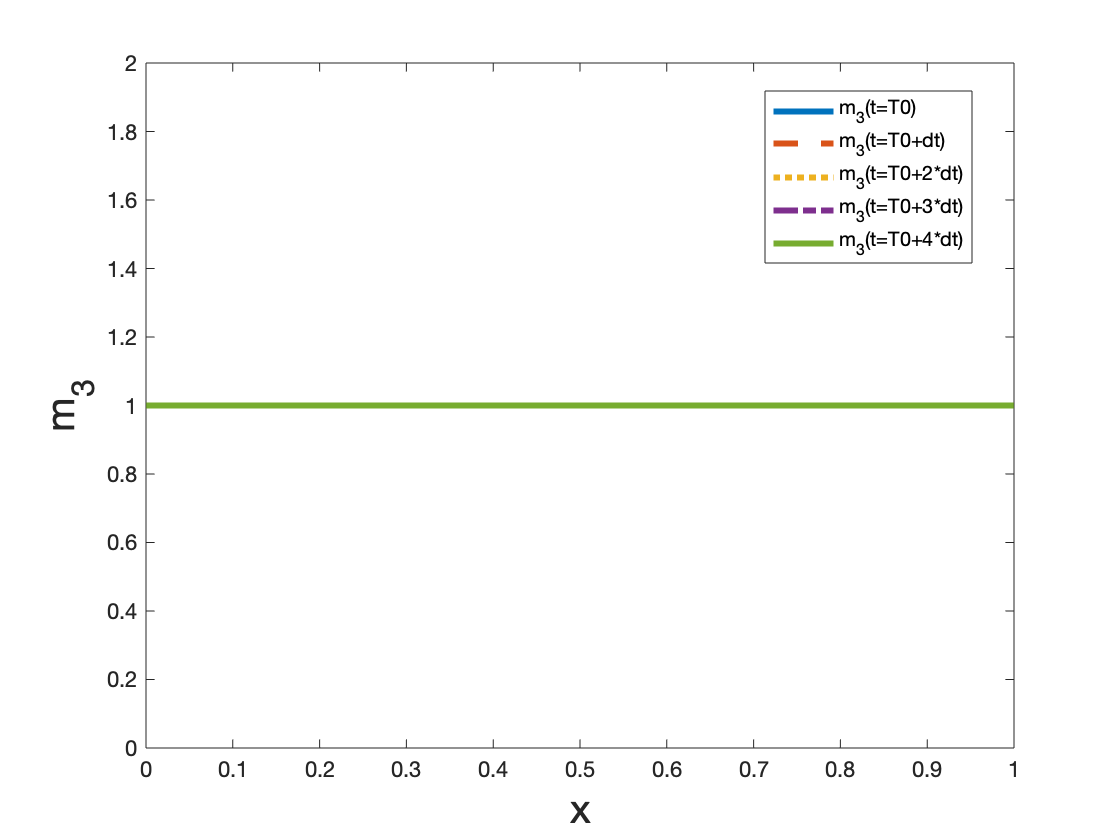}}
    \subfloat[$T0=0.01$, $m_3$]{\includegraphics[width=0.5\linewidth]{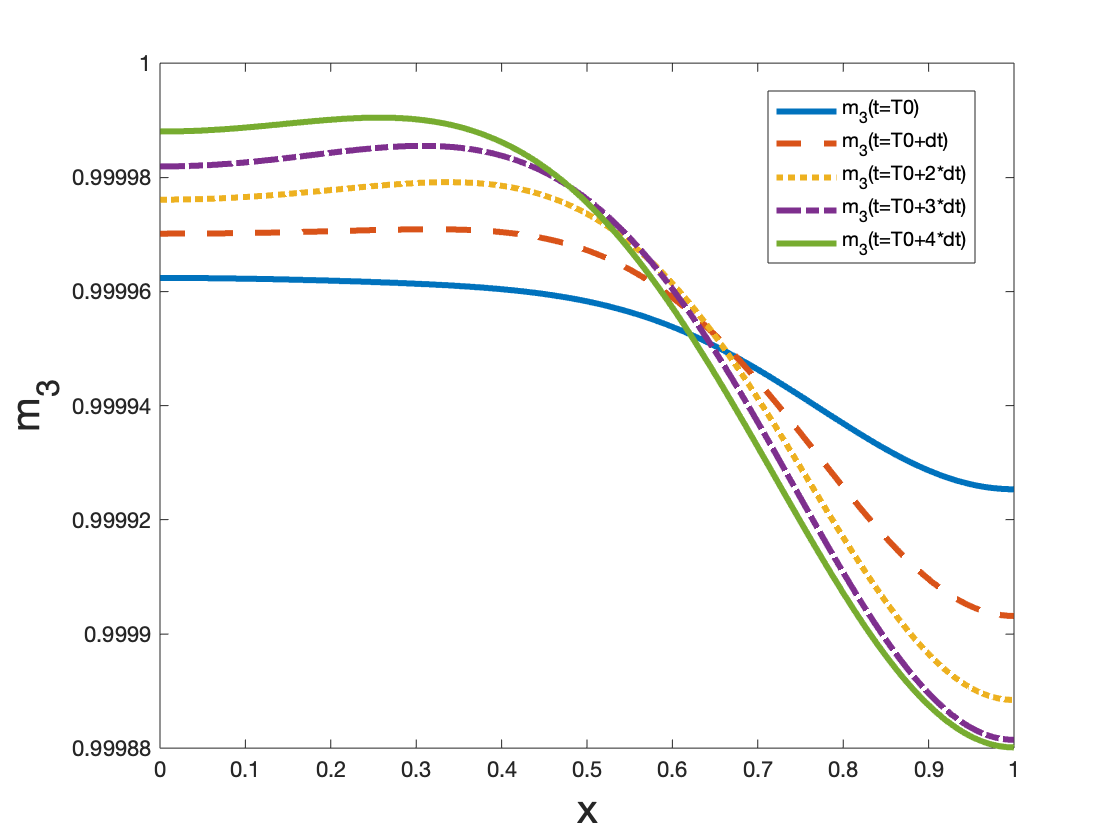}}
    \caption{The solution profile using proposed method in 1D given the initial condition $m_0$ with $T0$ specified without source term, $\alpha=0.01$ and $T=0.1$, $N_x=2000$, $N_t=5$.}
    \label{fig:1}
\end{figure}

\begin{figure}[htbp]
    \centering
    \subfloat[$T0=0$, $m_1$]{\includegraphics[width=0.5\linewidth]{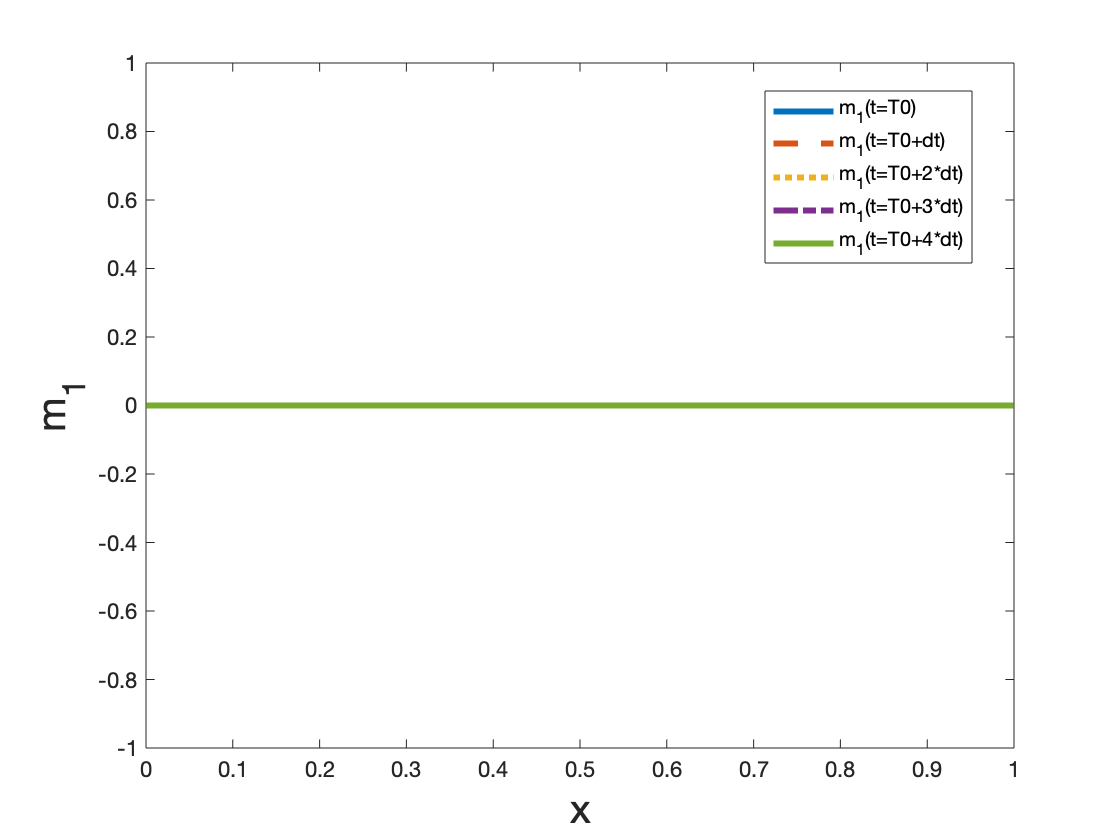}}
    \subfloat[$T0=0.01$, $m_1$]{\includegraphics[width=0.5\linewidth]{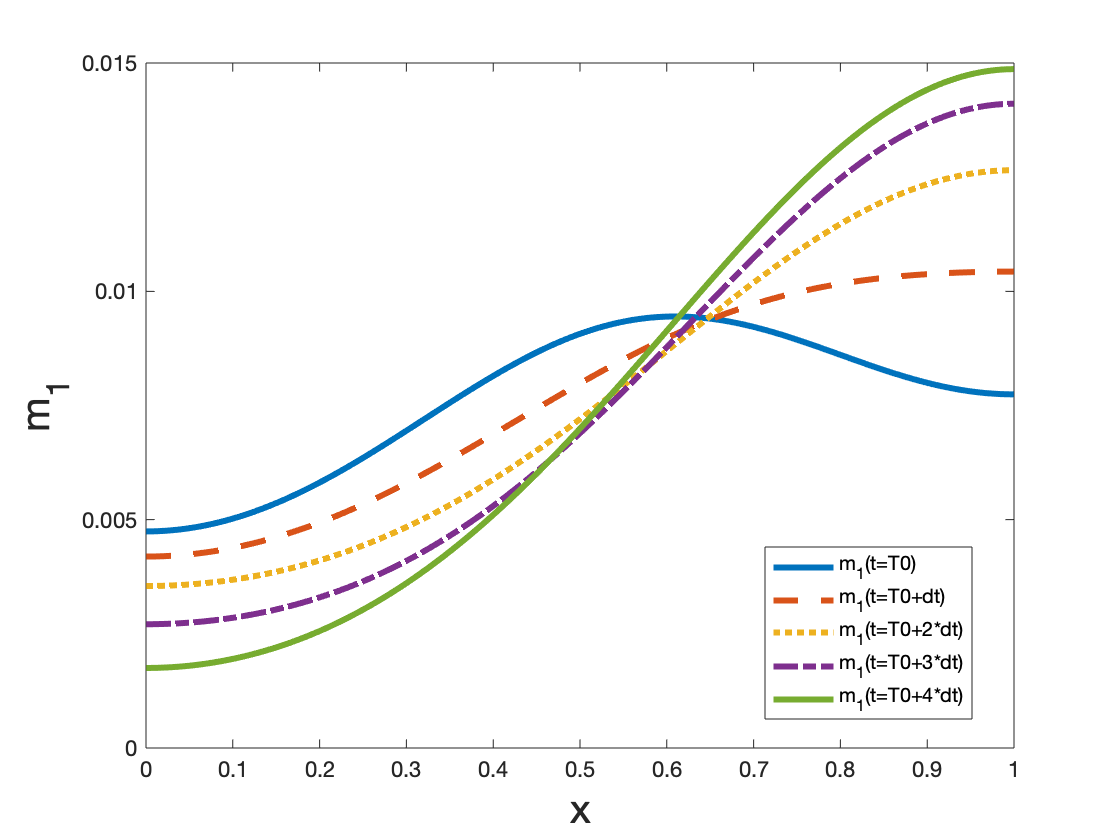}}
    \hspace{0.1in}
    \subfloat[$T0=0$, $m_2$]{\includegraphics[width=0.5\linewidth]{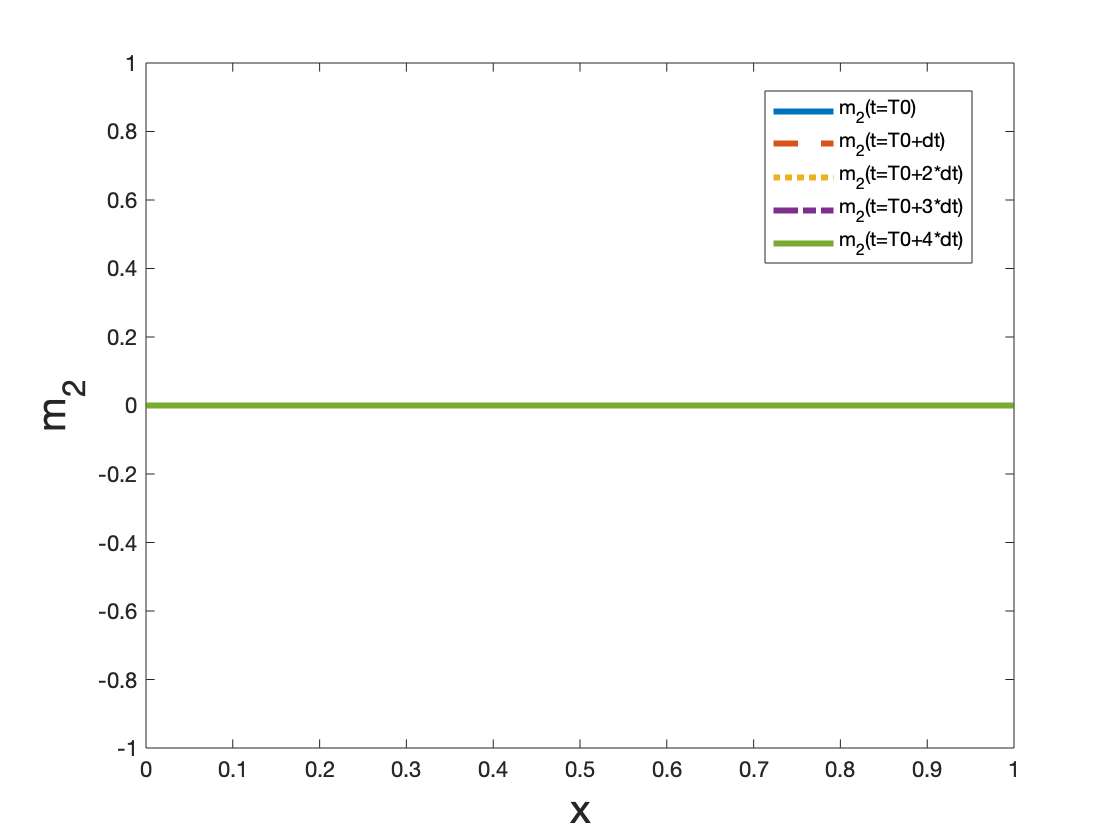}}
    \subfloat[$T0=0.01$, $m_2$]{\includegraphics[width=0.5\linewidth]{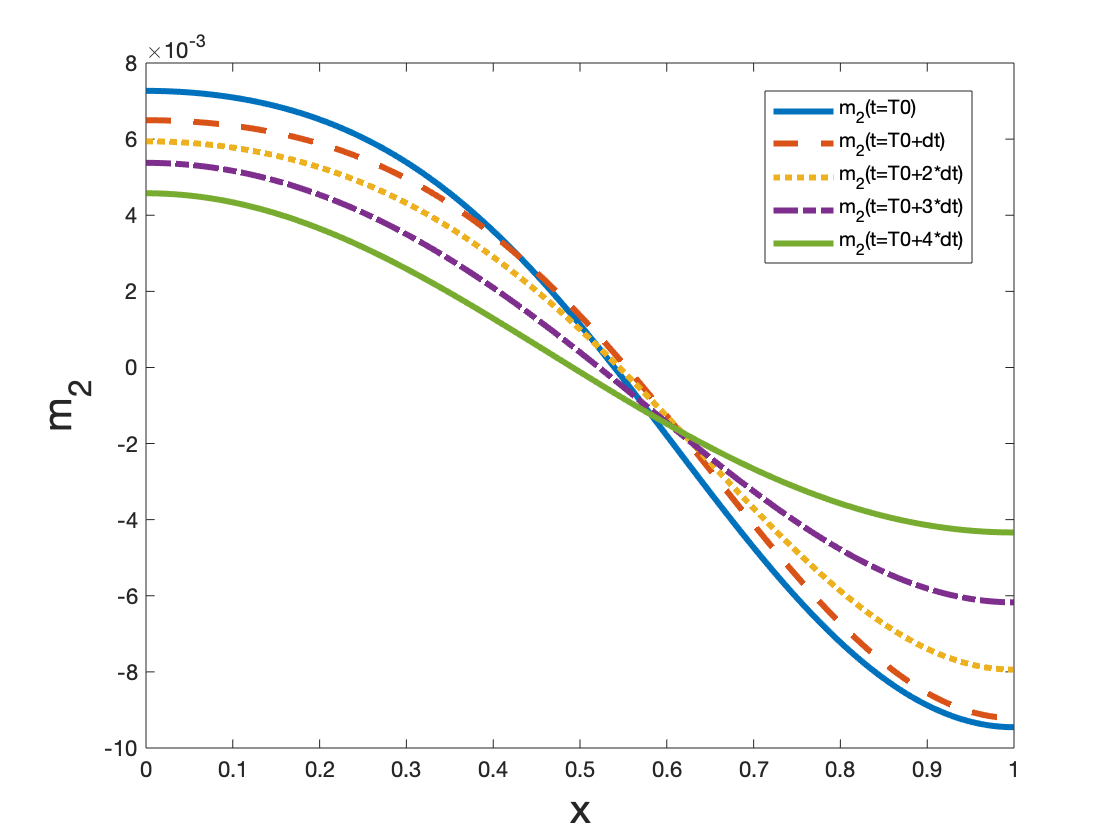}}
    \hspace{0.1in}
    \subfloat[$T0=0$, $m_3$]{\includegraphics[width=0.5\linewidth]{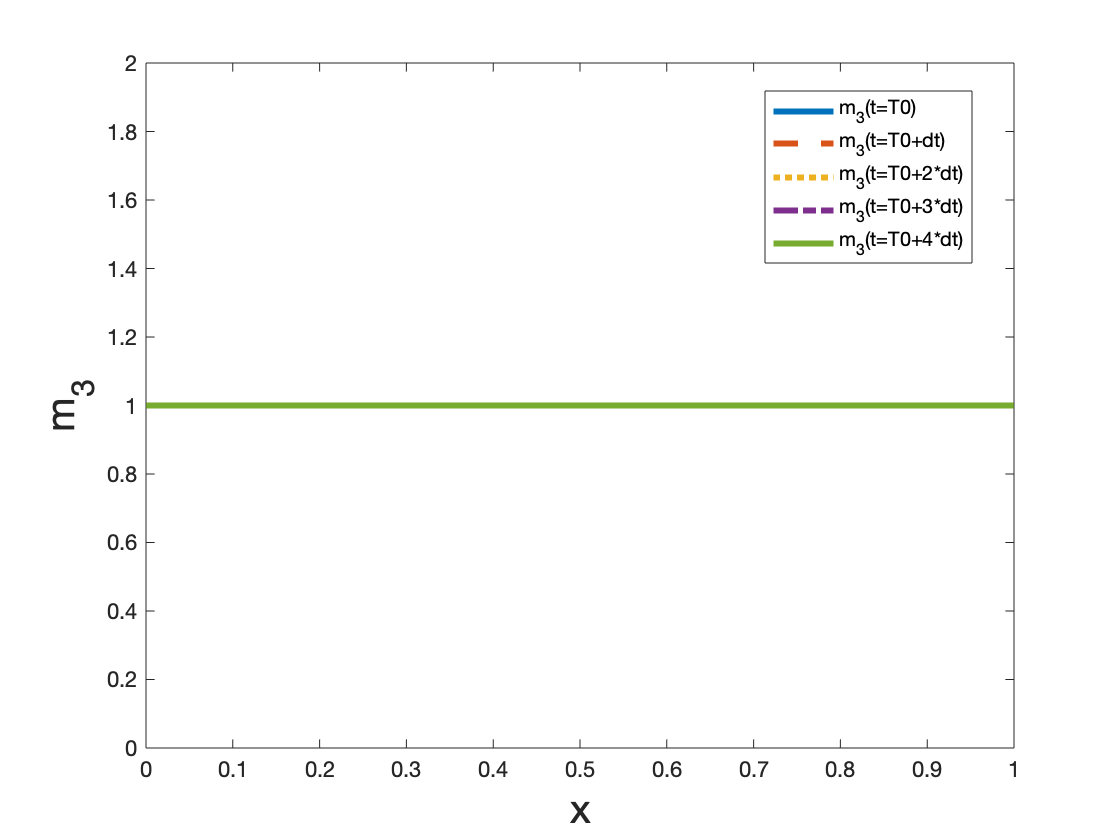}}
    \subfloat[$T0=0.01$, $m_3$]{\includegraphics[width=0.5\linewidth]{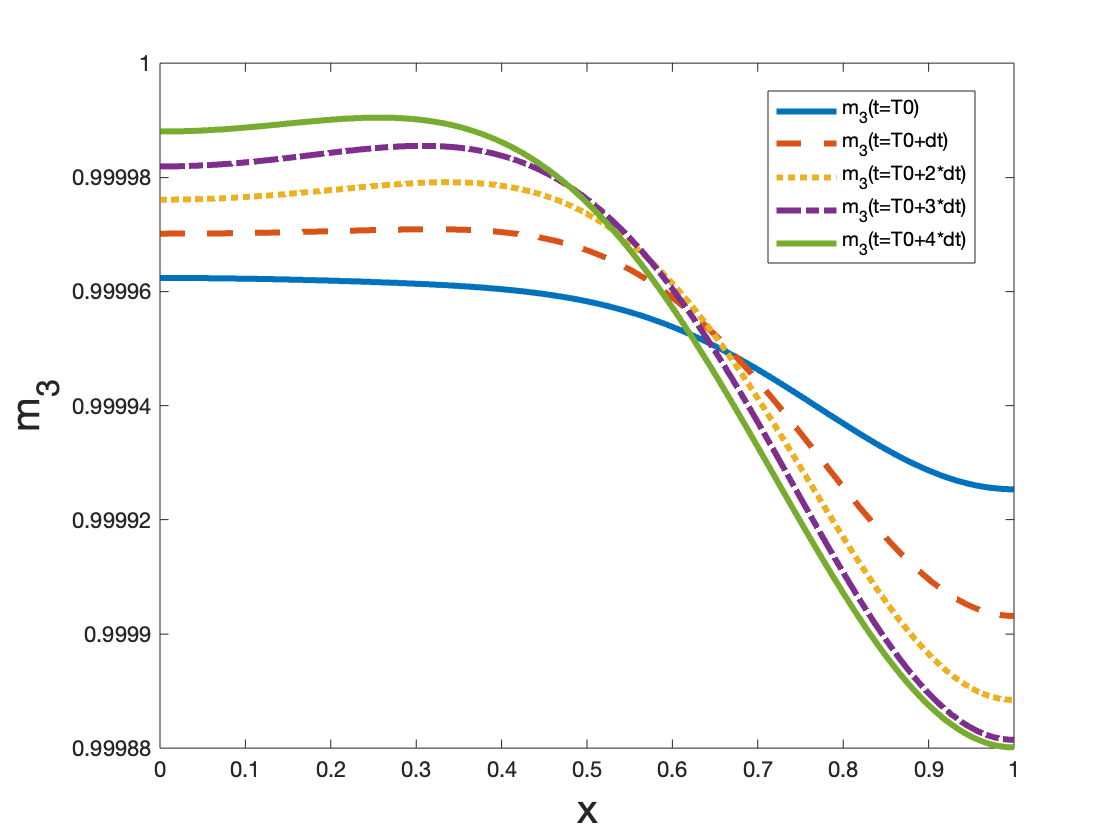}}
    \caption{The solution profile using BDF1 projection method in 1D given the initial condition $m_0$ with $T0$ specified without source term, $\alpha=0.01$ and $T=0.1$, $N_x=2000$, $N_t=5$.}
    \label{fig:2}
\end{figure}

To demonstrate the numerical consistency of the proposed method, comparative simulations are performed using both the proposed scheme and the first-order BDF projection method. The corresponding results for the 1D case are presented in \Cref{fig:1} and \Cref{fig:2}, with the simulation parameters set as \(\alpha=0.01\), \(N_x=2000\), and \(N_t=5\).

For the 3D scenario, the initial conditions selected to verify numerical consistency are:

\begin{align*}
    \m_0&=\left(\cos(\cos(\pi x))\sin (0.01), \sin(\cos(\pi x))\sin (0.01), \cos (0.01)\right)^T\\
    \m_0&=\left(\cos(\cos(\cos(\pi x)))\sin (\pi x+t), \sin(\cos(\cos(\pi x)))\sin (\pi x+t), \cos (\pi x+t)\right)^T
\end{align*}

\begin{figure}[htbp]
    \centering
    \subfloat[arrow profile]{\includegraphics[width=0.5\linewidth]{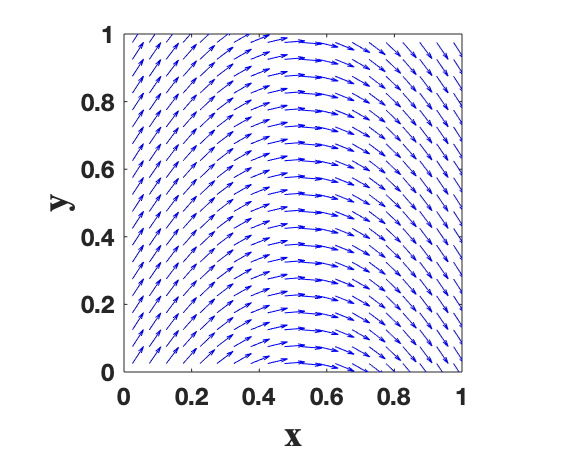}}
    \subfloat[angle profile]{\includegraphics[width=0.53\linewidth]{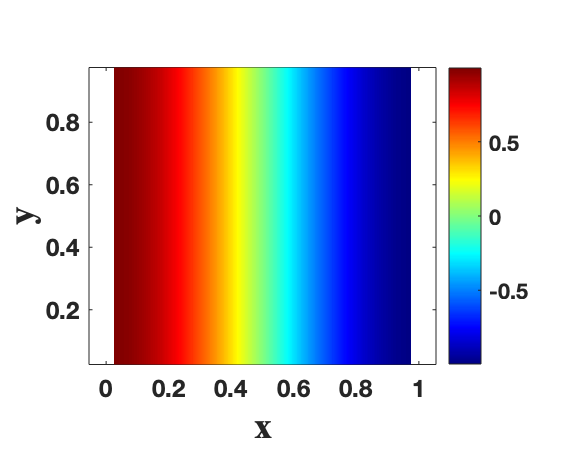}}
    \hspace{0.1in}
    \subfloat[arrow profile]{\includegraphics[width=0.5\linewidth]{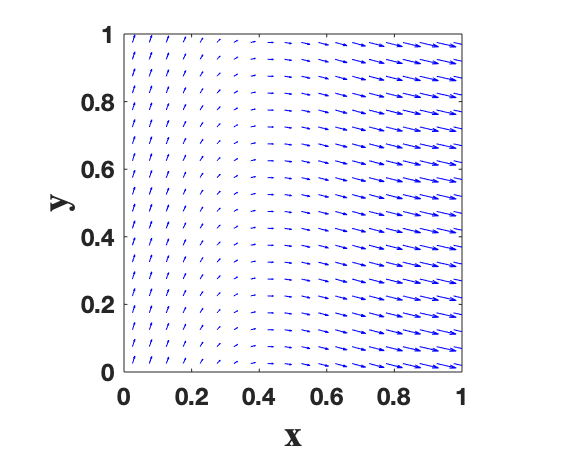}}
    \subfloat[angle profile]{\includegraphics[width=0.53\linewidth]{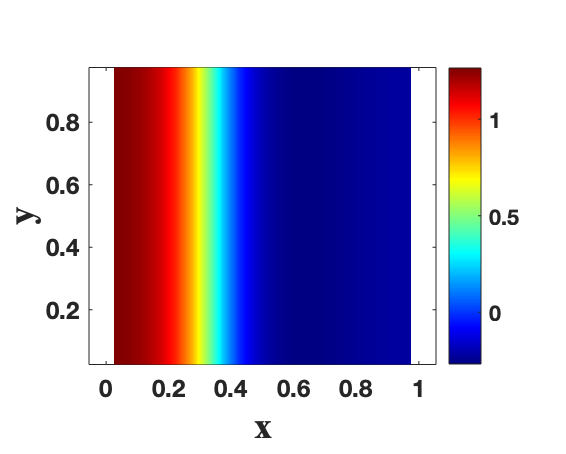}}
    \hspace{0.1in}
    \subfloat[arrow profile]{\includegraphics[width=0.5\linewidth]{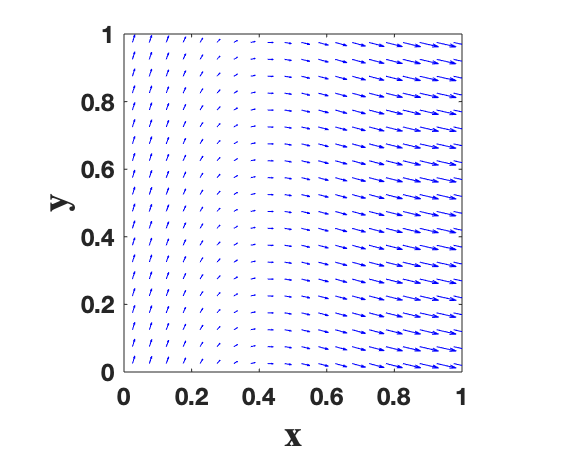}}
    \subfloat[angle profile]{\includegraphics[width=0.53\linewidth]{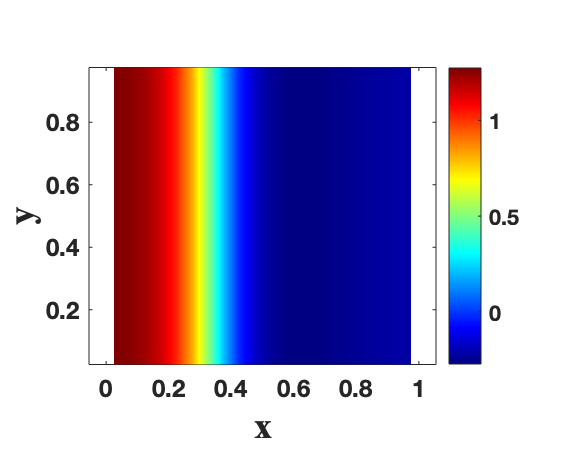}}
      \caption{The solution profile using GSPM and proposed methods in 3D given the initial condition $m_0$ with initial condition specified without source term, $\alpha=0$ and $T=0.1$, $N_x=N_y=N_z=20$, $N_t=400$. Top row with initial condition; Middle row with GSPM; Bottom row with proposed method. Initial condition given: $\m_0=[\cos(\cos(\pi x))\sin(0.01),\sin(\cos(\pi x))\sin(0.01),\cos(0.01)]$}
    \label{fig:3}
\end{figure}

\begin{figure}[htbp]
    \centering
    \subfloat[arrow profile]{\includegraphics[width=0.5\linewidth]{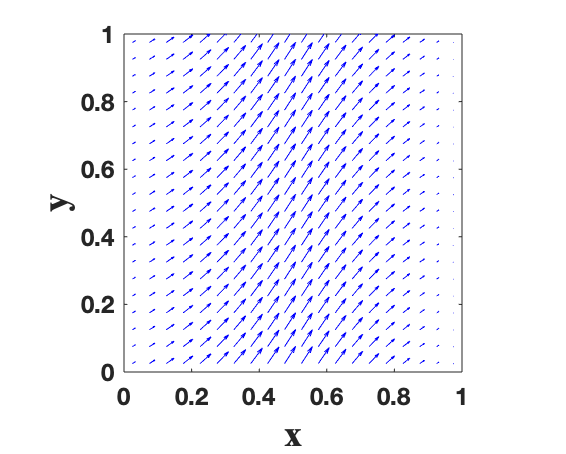}}
    \subfloat[angle profile]{\includegraphics[width=0.53\linewidth]{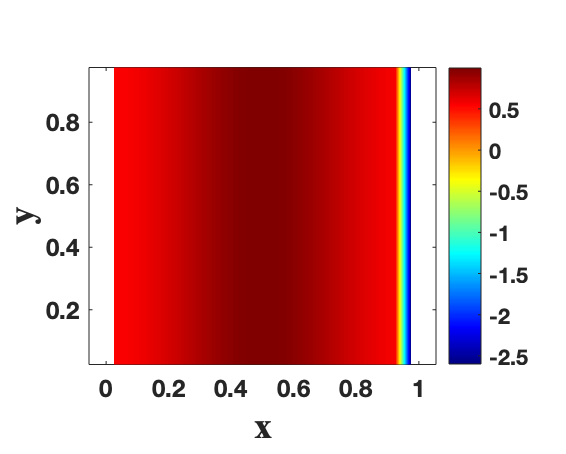}}
    \hspace{0.1in}
    \subfloat[arrow profile]{\includegraphics[width=0.5\linewidth]{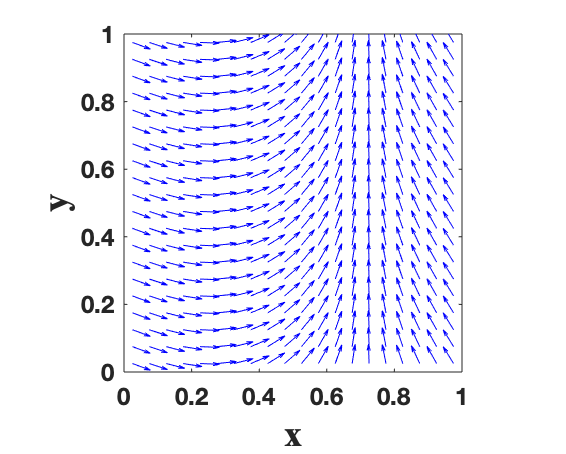}}
    \subfloat[angle profile]{\includegraphics[width=0.53\linewidth]{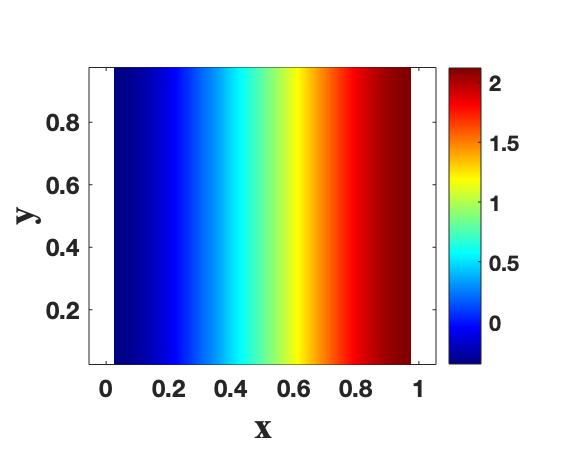}}
    \hspace{0.1in}
    \subfloat[arrow profile]{\includegraphics[width=0.5\linewidth]{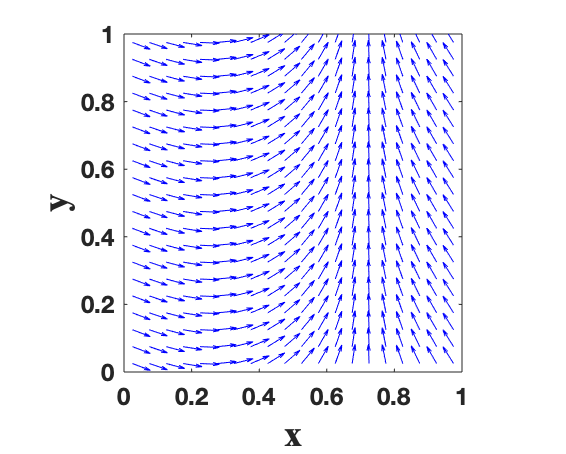}}
    \subfloat[angle profile]{\includegraphics[width=0.53\linewidth]{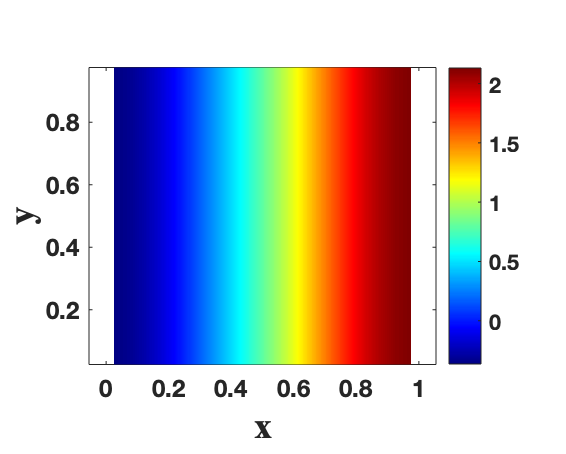}}
      \caption{The solution profile using GSPM and proposed methods in 3D given the initial condition $m_0$ with initial condition specified without source term, $\alpha=0$ and $T=0.1$, $N_x=N_y=N_z=20$, $N_t=40$. Top row with initial condition; Middle row with BDF1; Bottom row with proposed method. Initial condition given: $\m_0=[\cos(\cos(\cos(\pi x)))\sin(\pi x+t),\sin(\cos(\cos(\pi x)))\sin(\pi x+t),\cos(\pi x+t)]$ with $t=T0=0$.}
    \label{fig:6}
\end{figure}

The numerical results corresponding to these 3D initial conditions are presented in \Cref{fig:3} and \Cref{fig:6}, which effectively verify the numerical consistency of the proposed method.

To assess the robustness of the proposed method under diverse initial configurations, additional initial conditions are tested, which are given by:

\begin{align*}
        \m_0&=\left(\cos(XYZ)\sin (0.01), \sin(XYZ)\sin (0.01), \cos (0.01)\right)^T\\
         \m_0&=\left(\cos(\cos(\pi x)\cos(\pi y)\cos(\pi z))\sin (0.01), \sin(\cos(\pi x)\cos(\pi y)\cos(\pi z))\sin (0.01), \cos (0.01)\right)^T.
\end{align*}

\begin{figure}[htbp]
    \centering
    \subfloat[arrow profile]{\includegraphics[width=0.5\linewidth]{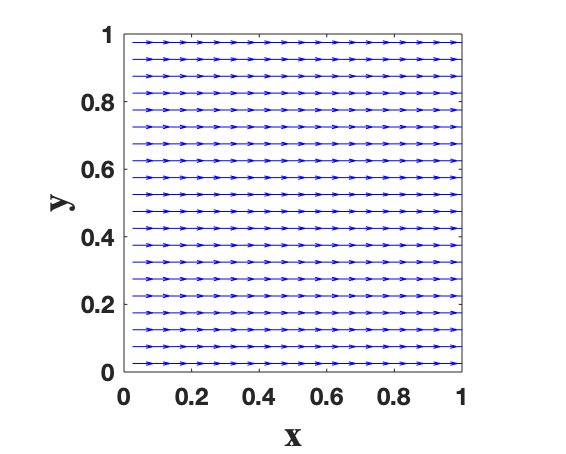}}
    \subfloat[angle profile]{\includegraphics[width=0.53\linewidth]{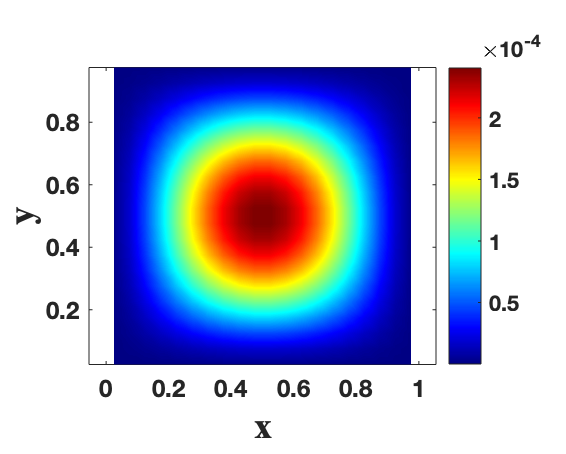}}
    \hspace{0.1in}
    \subfloat[arrow profile]{\includegraphics[width=0.5\linewidth]{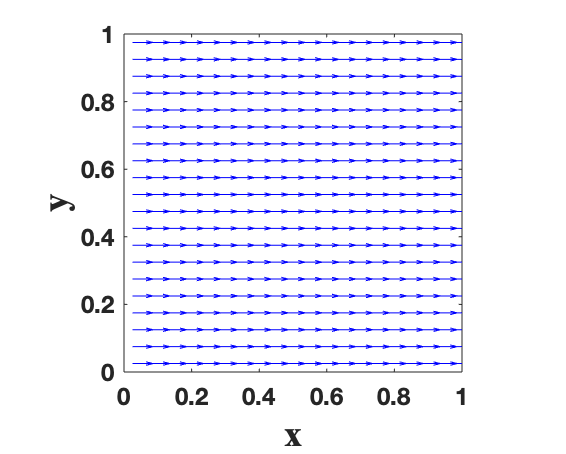}}
    \subfloat[angle profile]{\includegraphics[width=0.53\linewidth]
    {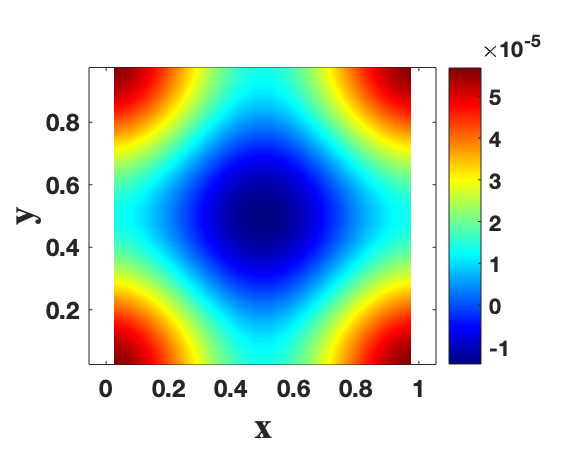}}
    \hspace{0.1in}
     \subfloat[arrow profile]{\includegraphics[width=0.5\linewidth]{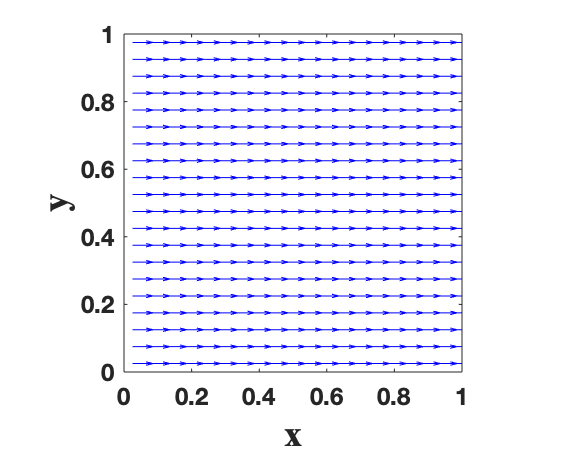}}
    \subfloat[angle profile]{\includegraphics[width=0.53\linewidth]{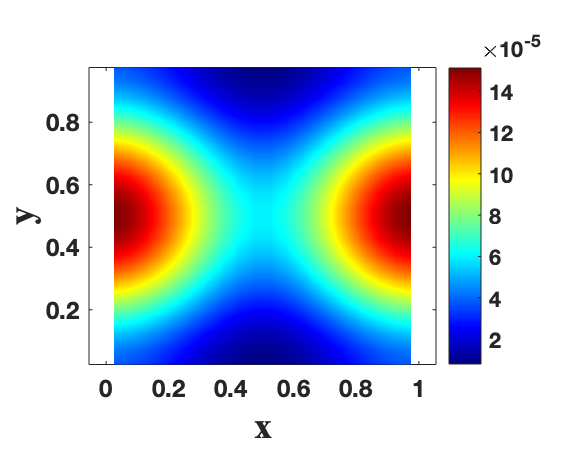}}
      \caption{Given initial conditions in 3D as $\m_0=[\cos(XYZ)\sin(0.01),\sin(XYZ)\sin(0.01),\cos(0.01)]$ specified without source term, $\alpha=0.01$ and $T=0.1$, $N_x=N_y=N_z=20$, $N_t=40$. Top row: initial condition; Middle row: BDF1 projection method; Bottom row: proposed method.}
    \label{fig:4}
\end{figure}

\begin{figure}[htbp]
    \centering
    \subfloat[arrow profile]{\includegraphics[width=0.5\linewidth]{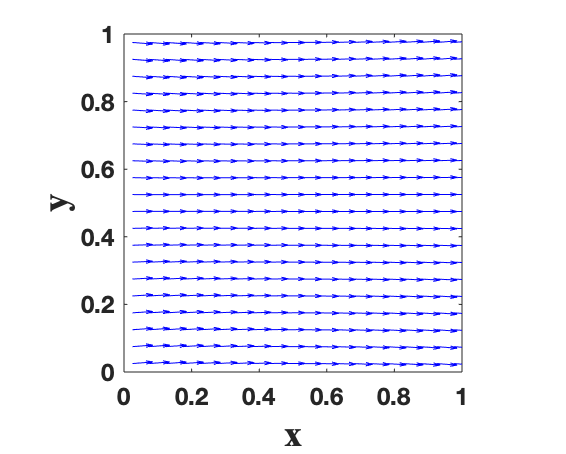}}
    \subfloat[angle profile]{\includegraphics[width=0.53\linewidth]{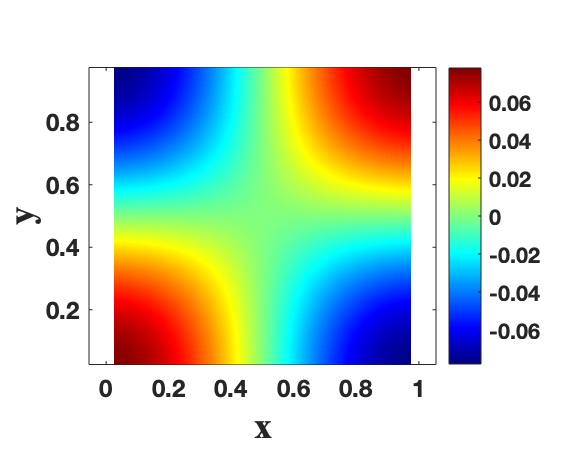}}
     \hspace{0.1in}
    \subfloat[arrow profile]{\includegraphics[width=0.5\linewidth]{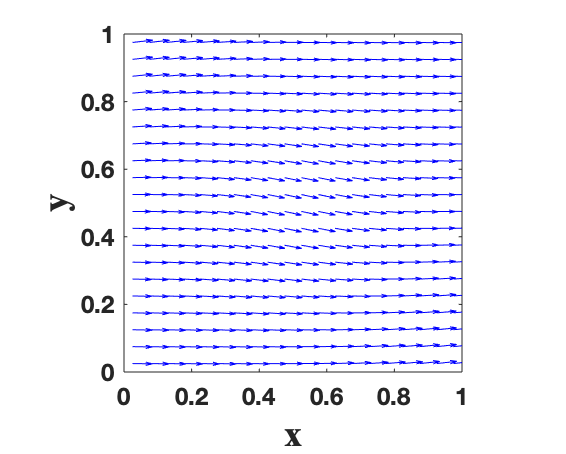}}
    \subfloat[angle profile]{\includegraphics[width=0.53\linewidth]{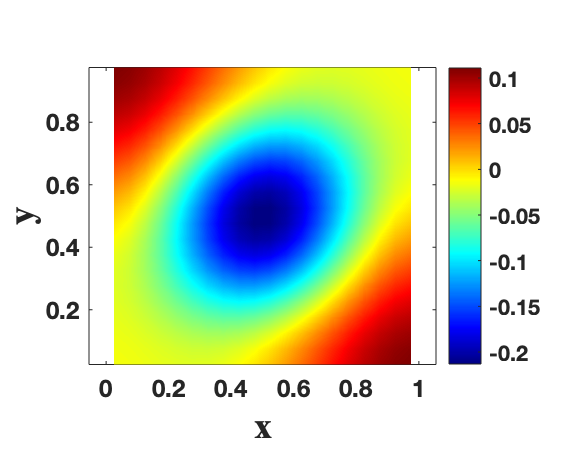}}
     \hspace{0.1in}
    \subfloat[arrow profile]{\includegraphics[width=0.5\linewidth]{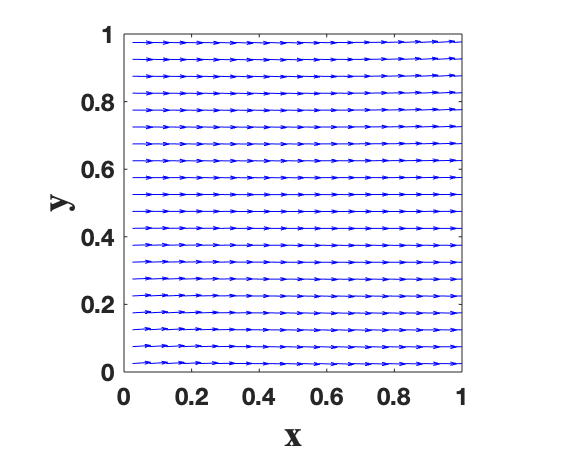}}
    \subfloat[angle profile]{\includegraphics[width=0.53\linewidth]{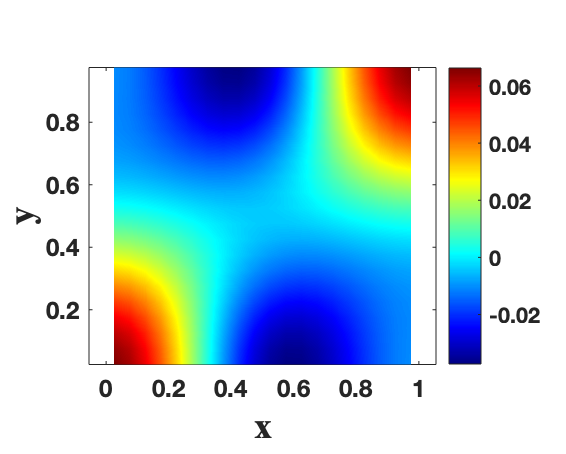}}
      \caption{Given initial conditions in 3D as $\m_0=[\cos(\cos(\pi x)\cos(\pi y)\cos(\pi z))\sin(0.01),\sin(\cos(\pi x)\cos(\pi y)\cos(\pi z))\sin(0.01),\cos(0.01)]$ specified without source term, $\alpha=0.01$ and $T=0.1$, $N_x=N_y=N_z=20$, $N_t=40$. Top row: initial condition; Middle row: BDF1 projection method; Bottom row: proposed method.}
    \label{fig:5}
\end{figure}

Here, \(X=x^2(1-x)^2\), \(Y=y^2(1-y)^2\), and \(Z=z^2(1-z)^2\) are the same spatial modulation terms as defined in the 3D exact solution. The simulation results for these additional initial conditions are presented in \Cref{fig:4} and \Cref{fig:5}, confirming both the stability and robustness of the proposed method.

\section{Conclusions and discussions}
\label{sec:conclusions}

In this paper, a structure-preserving numerical method is proposed for solving the Landau-Lifshitz-Gilbert (LLG) equation, which achieves first-order accuracy in time and second-order accuracy in space. This method integrates two key iterative procedures: a first-order semi-implicit Backward Differentiation Formula (BDF) iteration and a Crank-Nicolson-type iteration. Specifically, the first step employs the semi-implicit BDF iteration to provide a stable configuration for the higher-order terms in the LLG equation, while the second step utilizes the Crank-Nicolson-type iteration to preserve the norm constraint of the magnetization field. A distinct advantage of the proposed method is that it is constructed solely based on the intrinsic structure of the LLG equation, eliminating the need for an additional projection step—an approach that avoids the nonlinearity-induced challenges in stability and convergence analysis associated with conventional projection-based methods.
Several directions for future work are outlined as follows. First, we will conduct a rigorous theoretical analysis of the stability of the proposed method and further validate its effectiveness through applications in practical micromagnetic simulations. Second, we will modify the proposed method to adopt the finite element method for spatial discretization, enhancing its adaptability to complex geometric domains. Third, the method can be extended to higher-order accuracy by leveraging BDF-k schemes with \(k=2,3,4,5\). Additionally, the convergence analysis and stability analysis of the proposed method, together with the associated normalizing step, will be rigorously proven in subsequent research. Overall, the proposed method exhibits excellent performance in terms of numerical stability, magnetization length preservation, and computational accuracy, laying a solid foundation for its broader application in micromagnetics.

\section*{Acknowledgments}
This work is supported in part by the Jiangsu Science and Technology Programme-Fundamental Research Plan Fund (BK20250468), and the Research and Development Fund of Xi'an Jiaotong Liverpool University (RDF-24-01-015).

\vspace{1cm}

\bibliographystyle{elsarticle-num-names}
\bibliography{references.bib}

\end{document}